\documentclass[12pt,fleqn,epsfig]{article}

\usepackage{graphicx}
\usepackage{amssymb}

\newtheorem{theorem}{Theorem}[section]
\newtheorem{proposition}[theorem]{Proposition}
\newtheorem{lemma}[theorem]{Lemma}

\newtheorem{claim}{Claim}[theorem]
\newtheorem{conjecture}[theorem]{Conjecture}

\begin{document}

\title{Intersecting 1-factors and nowhere-zero 5-flows}

\author{Eckhard Steffen\thanks{Paderborn Institute for Advanced Studies in 
		Computer Science and Engineering,
		Paderborn University,	
		Zukunftsmeile 1, 33102 Paderborn, Germany; 
		es@upb.de}}

\date{}

\maketitle

\begin{abstract}
{\small{Let $G$ be a bridgeless cubic graph, and $\mu_2(G)$ the minimum number $k$ such that
two 1-factors of $G$ intersect in $k$ edges. A cyclically $n$-edge-connected cubic graph $G$ has a nowhere-zero 5-flow if 
(1) $n \geq 6$ and $\mu_2(G) \leq 2$ or (2) if $n \geq 5 \mu_2(G)-3$.
}}
\end{abstract}  

\section[]{Introduction}

One of the first and famous
theorems of graph theory, Petersen's Theorem from 1891, states that every bridgeless cubic graph has a 1-factor, and hence 
a 2-factor as well. The edge-chromatic number $\chi'(G)$ of a cubic graph $G$ is either 3 or 4.
 If a cubic graph has two disjoint 1-factors, then it is 3-edge-colorable. 
Hence, if $\chi'(G)=4$, then any two 1-factors of $G$ have a non-empty intersection.  
Let $\mu_2(G) = \min \{|M_1 \cap M_2| : M_1 \mbox{ and } M_2 \mbox{ are 1-factors of } G\}$. 

An integer nowhere-zero $k$-flow on a graph $G$ is an assignment of a direction and a value of $\{1, \dots, (k-1)\}$ to each edge of
$G$ such that the Kirchhoff's law is satisfied at every vertex of $G$.  (This is the most restrictive definition of a nowhere-zero $k$-flow.
But it is equivalent to more flexible definitions, see e.g.~\cite{Seymour_95}.)
A cubic graph $G$ is bipartite if and only if it has a nowhere-zero 3-flow, and $\chi'(G)=3$  if and only if $G$ has 
a nowhere-zero 4-flow. Seymour \cite{Seymour_81} proved that every bridgeless graph has a nowhere-zero 6-flow. 
So far this is the best approximation to Tutte's famous 5-flow conjecture, which is equivalent to its restriction to cubic graphs.

\begin{conjecture} [\cite{Tutte_54}] \label{5FC}
Every bridgeless graph has a nowhere-zero 5-flow.
\end{conjecture}

Kochol \cite{Kochol_04} proved that a minimum counterexample to the 5-flow conjecture is a cyclically 6-edge-connected cubic graph. Hence it suffices to 
prove Conjecture \ref{5FC} for these graphs. The following theorem is the main result of this paper. 

\begin{theorem} \label{gamma2_5flow}
Let $G$ be a cyclically 6-edge-connected cubic graph. If $\mu_2(G) \leq 2$, then $G$ has a nowhere-zero 5-flow.
\end{theorem}

In section \ref{gamma_2} we observe that $2\mu_2(G) \geq \omega(G)$, where $\omega(G)$ is the oddness of $G$, which is the 
minimum number of odd circuits in a 2-factor of $G$. Jaeger \cite{Jaeger_88} showed that cubic graphs with oddness at most 2 have a nowhere-zero 5-flow.
Using results of \cite{Steffen_10} we further deduce the following theorem.

\begin{theorem} \label{gamma_cc}
Let $G$ be a cyclically $k$-edge-connected cubic graph. If $k \geq 5 \mu_2(G) - 3$,
then $G$ has a nowhere-zero 5-flow.  
\end{theorem}

Theorems \ref{gamma2_5flow} and \ref{gamma_cc} are consequences of the stronger Theorems \ref{gamma2_5flow_sharp} and 
\ref{5flow_cc}, respectively.

\section[] {Nowhere-zero 5-flows} \label{gamma_2}

Let $G$ be a graph and $S \subseteq V(G)$. The set of edges with precisely one end in $S$ is
denoted by $\partial_G(S)$. We start with the following folklore result.

\begin{proposition} \label{min_intersecting_oddness}
If $G$ is a bridgeless cubic graph, then $\omega(G) \leq 2 \mu_2(G)$.
\end{proposition}
{\bf Proof.} Let $\mu_2(G) \geq 0$, and $M_1$, $M_2$ be two 1-factors of $G$ with $|M_1 \cap M_2|= \mu_2(G)$.
The 2-factor $\overline{M_1}$, which is complementary to $M_1$, has at least $\omega(G)$ odd circuits.
If $C$ is an odd circuit of $\overline{M_1}$, then $\partial_G(V(C)) \subset M_1$. Since $|\partial_G(V(C))|$ is odd
it follows that $\partial_G(V(C)) \cap M_2 \not = \emptyset$. 
Hence $\mu_2(G) \geq \omega(G)/2$.  \hfill $square$

A minimum 2-factor of a cubic graph $G$ has precisely $\omega(G)$ odd circuits.
Let $\omega(G) \geq 2$, ${\cal F}$ be a minimum 2-factor, and 
$m_G({\cal F})$ be the maximum number $k$ such 
that $G$ has no edge cut $E$ with fewer than $k$ edges such that two
components of $G - E$ contain odd circuits of ${\cal F}$. 
Let $m_G^* = \infty$ if $\omega(G) = 0$, and 
$m_G^* = \max \{ m_G({\cal F}) : {\cal F} \mbox{ is a minimum 2-factor of } G \}$
if $\omega(G) > 0$. 
 
\begin{theorem} [\cite{Steffen_10}]\label{Main}
Let $G$ be a bridgeless cubic graph. If $m_G^* \geq \frac{5}{2} \omega(G) - 3$,
then $G$ has a nowhere-zero 5-flow.
\end{theorem} 

Thus the following theorem follows with Proposition \ref{min_intersecting_oddness}.

\begin{theorem} \label{5flow_cc} Let $G$ be a bridgeless cubic graph.
If $m_G^* \geq 5 \mu_2(G) - 3$, then $G$ has a nowhere-zero 5-flow. 
\end{theorem}

Since the cyclic connectivity of a cubic graph is smaller or equal to $m_G^*$, Theorem \ref{5flow_cc} implies 
Theorem \ref{gamma_cc}.

We now consider graphs where $\mu_2(G)$ is small. 
An {\em orientation} $D$ of $G$ is an assignment of a
direction to each edge. For $S \subseteq V(G)$, 
$D^-(S)$ ($D^+(S)$) is the set of edges of $\partial_G(S)$ whose head 
(tail) is incident to a vertex of $S$.
The oriented graph is denoted by $D(G)$, 
$d_{D(G)}^-(v) = |D^-(\{v\})|$ and $d_{D(G)}^+(v) = |D^+(\{v\})|$ denote the {\em indegree}
and {\em outdegree} of vertex $v$ in $D(G)$, respectively. The 
degree of a vertex $v$ in the undirected graph $G$ is $d_{D(G)}^+(v) + d_{D(G)}^-(v)$, and it is denoted by $d_G(v)$.  

Let $k$ be a positive integer, and $\varphi$  a function from the edge set of the directed graph $D(G)$
into the set $\{0, 1, \dots, (k-1)\}$. For $S \subseteq V(G)$ let 
$\delta \varphi (S) = \sum_{e \in D^+(S)}\varphi(e) - \sum_{e \in D^-(S)}\varphi(e)$. 
The function $\varphi$ is a $k$-flow on $G$ if $\delta \varphi(S) = 0$ for every $S \subseteq V(G)$. 
 The {\em support} of $\varphi$  is 
the set $\{e \in E(G) : \varphi(e) \not = 0\}$, and it is denoted by $supp(\varphi)$.
A $k$-flow $\varphi$ is a nowhere-zero $k$-flow if $supp(\varphi) = E(G)$.

We will use balanced valuations of graphs, which were introduced by 
Bondy \cite{Bondy} and Jaeger \cite{Jaeger_75}. 
A {\em balanced valuation} of a graph
$G$ is a function $f$ from the vertex set $V(G)$ into the real numbers, such that
$| \sum_{v \in X} f(v) | \leq | \partial_G(X) |$ for all $X \subseteq V(G)$. We will use
the following fundamental theorem of Jaeger.

\begin{theorem} [\cite{Jaeger_75}] \label{Thm_Jaeger_75} 
Let $M$ be a graph with orientation $D$ and $k\geq 3$. Then $M$ 
has a nowhere-zero $k$-flow if and only if there 
is a balanced valuation $f$ of $M$ with 
$ f(v) = \frac{k}{k-2}(2d_{D(M)}^+(v) - d_M(v))$, for all $v \in V(M).$
\end{theorem}

In particular, Theorem \ref{Thm_Jaeger_75} says that a cubic graph $G$ has a nowhere-zero 4-flow 
(nowhere-zero 5-flow) if and only if there is a balanced valuation of $G$ with values in 
$\{ \pm 2\}$ ($\{ \pm \frac{5}{3}\}$).

For the following we define a specific nowhere-zero 4-flow and the corresponding balanced valuation which are induced 
by a proper 3-edge-coloring of $G$. If we describe a flow which relies on a specific orientation $D$ of the edges of $G$, 
then we also write $(D,\varphi)$ instead of $\varphi$.
For $i \in\{1,2\}$ let $(D_i,\varphi_i)$ be flows on $G$.
The sum $(D_1,\varphi_1) + (D_2,\varphi_2)$ is the flow $(D, \varphi)$
on $G$ with orientation $D = D_1|_{ \{e : \varphi_1(e)   \geq  \varphi_2(e)\}}    \cup 
          D_2|_{\{e : \varphi_2(e) > \varphi_1(e)\}}$, and with flow value

\[ \varphi(e) = \left\{ \begin{array} {l@{\mbox{ }}l}
              \varphi_1(e) + \varphi_2(e), & \mbox{if } e  
              \mbox{ received the same direction
                in } D_1 \mbox{ and } D_2\\
              |\varphi_1(e) - \varphi_2(e)|, & \mbox{ otherwise. }
                      \end{array} \right. \]

Let $G$ be cubic graph, which has a proper 3-edge-coloring $c$. We define a canonical nowhere-zero 4-flow on $G$  as 
follows: For $i,j \in \{1,2,3\}$ with $1 \leq i < j \leq 3$ let $H_{i,j}$ be the even cycle which is induced by $c^{-1}(i) \cup c^{-1}(j)$.
Let $\phi_{1,2}$ be the flow on the directed circuits of $H_{1,2}$ with $\phi_{1,2}(e) = 1$ for all $e \in E (H_{1,2})$, and 
$\phi_{2,3}$ be the flow on the directed circuits of $H_{2,3}$ with $\phi_{2,3}(e) = 2$ for all $e \in E (H_{2,3})$. Then
$\phi = \phi_{1,2} + \phi_{2,3}$ is a nowhere-zero 4-flow on $G$. 
Note that the edges of $c^{-1}(1)$ have flow value 1, the edges of $c^{-1}(2)$ have flow value 1 or 3, and 
the edges of $c^{-1}(3)$ have flow value 2.

By Theorem \ref{Thm_Jaeger_75}, there is a balanced valuation $f$ of $G$ with $f(v) = 2 (2 d^+_{D(G)}(v) - d_{G}(v))$,
for $v \in V(G)$. It holds that $|2 d^+_{D(G)}(v) - d_{G}(v)| = 1$, and hence
$f(v) \in \{ \pm 2 \}$ for all $v \in V(G)$. 
The vertices of $G$ are
partitioned into two classes $A$ and $B$, where $A = \{ v : f(v) = 2\}$ and $B = \{v : f(v) = -2\}$. 
Call the elements of $A$ ($B$) the white (black) vertices of $G$. A balanced valuation which is induced by a  
canonical nowhere-zero 4-flow will be called a canonical balanced valuation of $G$. The following lemma
is an immediate consequence of the definition of a canonical nowhere-zero 4-flow.

\begin{lemma}  \label{verschiedene2}
Let $c$ be a proper 3-edge-coloring of a cubic graph $G$, and $A$, $B$ a partition of $V(G)$ which is induced by a canonical
nowhere-zero 4-flow with respect to $c$.  Let $x$, $y$ be the two vertices of an edge $e$. If $e \in c^{-1}(1) \cup c^{-1}(2)$, then 
$x$ and $y$ belong to different classes, i.e. $x \in A$ if and only if $y \in B$.
\end{lemma}

\begin{theorem} \label{gamma2_5flow_sharp}
Let $G$ be a cubic graph with $m_G^* \geq 6$. If $\mu_2(G) \leq 2$, then $G$ has a nowhere-zero 5-flow.
\end{theorem}
{\bf Proof.} If $\omega(G) \leq 2$, then the statement follows with Theorem \ref{5flow_cc}. Thus, by Proposition \ref{min_intersecting_oddness},
we have to consider the case when $\mu_2(G) = 2$ and $\omega(G) = 4$. It follows with Theorem \ref{5flow_cc} that $m_G^*(G) = 6$. 
Let $M_1$, $M_2$ be two 1-factors with minimum intersection and $M_1 \cap M_2 = \{e_1,e_2\}$, where $e_i = v_iw_i$ ($i \in \{1,2\}$). 
 
Let $G^-= G - \{e_1,e_2\}$ and $G'$ be the cubic graph which is obtained from $G^-$ by suppressing the bivalent vertices $v_i$, $w_i$. Let
$c'$ be a proper 3-edge-coloring of $G'$, such that the edges of $M_i$ are colored with color $i$.  

\begin{claim} \label{induced_2-factor}
There is a minimum 2-factor of $G$ with four odd circuits $C_{v_1}$, $C_{v_2}$,   $C_{w_1}$, $C_{w_2}$ such that
$x \in V(C_x)$ for $x \in \{v_1,v_2,w_1,w_2\}$.
\end{claim}
{\bf Proof.} Consider the 3-edge-coloring $c'$ of $G'$. Then each of $v_1$, $v_2$, $w_1$, $w_2$ subdivides an edge which is 
colored with color 3. Color one part of the subdivided edge with a new color 0, and $e_1$ and $e_2$ with color 1 to obtain a proper 4-edge-coloring 
$c$ of $G$. Then $G[c^{-1}(0) \cup c^{-1}(2) \cup c^{-1}(3)]$ is the desired 2-factor of $G$. Note that this 2-factor is the complement of $M_1$.
$\Box$

Let $\phi'$ be a canonical nowhere-zero 4-flow with respect to $c'$, and $A'$, $B'$ be the corresponding
partition of the vertices of $G'$. Subdividing an edge does not affect flow properties of graphs. Hence $\phi'$ induces 
a nowhere-zero 4-flow $\phi^-$ on $G^-$. 
Define $A = A' \cup \{v_1,v_2\}$, $B = B' \cup \{w_1,w_2\}$, and $f(v) = 5/3$, if 
$v \in A$ and $f(v) = - 5/3$ if $v \in B$. We claim, that $f$ is a balanced valuation of $G$ and hence $G$ has 
a nowhere-zero 5-flow by Theorem \ref{Thm_Jaeger_75}. 

Suppose to the contrary that $G$ does not have a nowhere-zero 5-flow; that is, $f$ is not a balanced valuation. 
Then there is a smallest $S \subset V(G)$ such that $|\sum_{v \in S}f(v)| > |\partial_G(S)|$.
Let $k$ be the difference between black and white vertices in $S$. Hence, $\frac{5}{3} k > |\partial_G(S)|$.

\begin{claim} \label{parity}
 $k \equiv |\partial_G(S)| (\bmod 2)$.
\end{claim} 
{\bf Proof.} If $k$ is even, then $|S \cap A|$ and $|S \cap B|$ have the same parity, and if $k$ is odd, then they have
different parities. Since $S$ is the disjoint union of $S \cap A$ and $S \cap B$  it follows that $k$ and $|S|$ have
the same parity. Since $G$ is cubic it follows that $k \equiv |\partial_G(S)| (\bmod 2)$. $\Box$

\begin{claim} \label{cut=6_k=4}
$|\partial_G(S)| = 6$ and $k=4$.
\end{claim} 
{\bf Proof.} If $|\{e_1, e_2\} \cap  \partial_G(S)| = 0$, then we obtain a contradiction, since -by the construction of $A$ and $B$- we
have $2k \leq |\partial_G(S)|$ in this case.

If $|\{e_1, e_2\} \cap  \partial_G(S)| = 1$, then, in the worst case, $2(k-1) \leq |\partial_G(S)| -1$ for the canonical balanced 
valuation of $G'$. Note that in this case $S=S'\cup\{x\}$, where $x \in \{v_1,v_2, w_1,w_2\}$, and $|\partial_G(S)| - 1 = |\partial_{G'}(S')|$.  
If $2(k-1) < |\partial_G(S)| -1$,
then $2k \leq |\partial_G(S)|$, a contradiction. Thus we may assume that $2(k-1) = |\partial_G(S)| -1$. Then $|\partial_G(S)|$ is odd,
and hence $k$ is odd by Claim \ref{parity}. If $|\partial_G(S)| \geq 5$, then 
$\frac{5}{3}k =  \frac{5}{6}2k = \frac{5}{6}|\partial_G(S)| + \frac{5}{6} \leq |\partial_G(S)|$, a contradiction.
Hence, $|\partial_G(S)| = 3$ which implies that $k=2$, a contradiction. 

If $|\{e_1, e_2\} \cap  \partial_G(S)| = 2$, then we similarly deduce a contradiction for the cases 
when  $2k$ or $2(k-1) \leq |\partial_G(S)|-2$. Hence it remains to consider the (worst) case, when $2(k-2) \leq |\partial_G(S)| - 2$
for the canonical balanced valuation of $G'$. 
If $|\partial_G(S)|$ is odd, then we obtain a contradiction as in the case above. Thus we can assume that $|\partial_G(S)|$ is even. 
If $2(k-2) < |\partial_G(S)| - 2$, then $2(k-2) \leq |\partial_G(S)| - 4$, and hence $2k \leq |\partial_G(S)|$, a contradiction.
It remains to consider the case when $2(k-2) = |\partial_G(S)| - 2$. Then $k = |\partial_G(S)|/2 + 1$. If $|\partial_G(S)| \geq 10$, then
$\frac{5}{3} k = \frac{5}{3} (|\partial_G(S)|/2 + 1) \leq |\partial_G(S)|$, a contradiction. 
Thus, $|\partial_G(S)| < 10$, and since $|\partial_G(S)|$ is even it follows with 
Claim \ref{parity} that $k$ is even. If $|\partial_G(S)| \in \{4,8\}$, then $k \in\{3,5\}$, a contradiction. Hence, $|\partial_G(S)| = 6$ and $k=4$.
$\Box$

Let $\partial_G(S) = \{e_1,e_2,f_1,f_2,f_3,f_4\}$, $f_i = x_iy_i$, and
$\{v_1, v_2, x_1, x_2, x_3, x_4 \} \subseteq S$. 
Let $\phi$ be the 4-flow on $G$ with $supp(\phi) = E(G) - \{e_1, e_2\}$, which is obtained from the (canonical) nowhere-zero 4-flow $\phi^-$ on $G^-$. 

Consider $G'$ and let $S'= S - \{v_1,v_2\}$. Since $k=4$, it follows that $k'= |S' \cap A'| - |S' \cap B'| = 2$. Lemma \ref{verschiedene2}
implies that $k' = |c'^{-1}(1) \cap \partial_{G'}(S')| = |c'^{-1}(2) \cap \partial_{G'}(S')|$.
Hence 
two edges of $\partial_G(S) - \{e_1,e_2\}$ are colored with color 1, say $f_1$, $f_2$ and two edges of $\partial_G(S) - \{e_1,e_2\}$ are colored with color 2, say $f_3$, $f_4$. The edges $e_1$ and $e_2$ are uncolored. 

\begin{figure}
\centering
 \includegraphics[width=4cm]{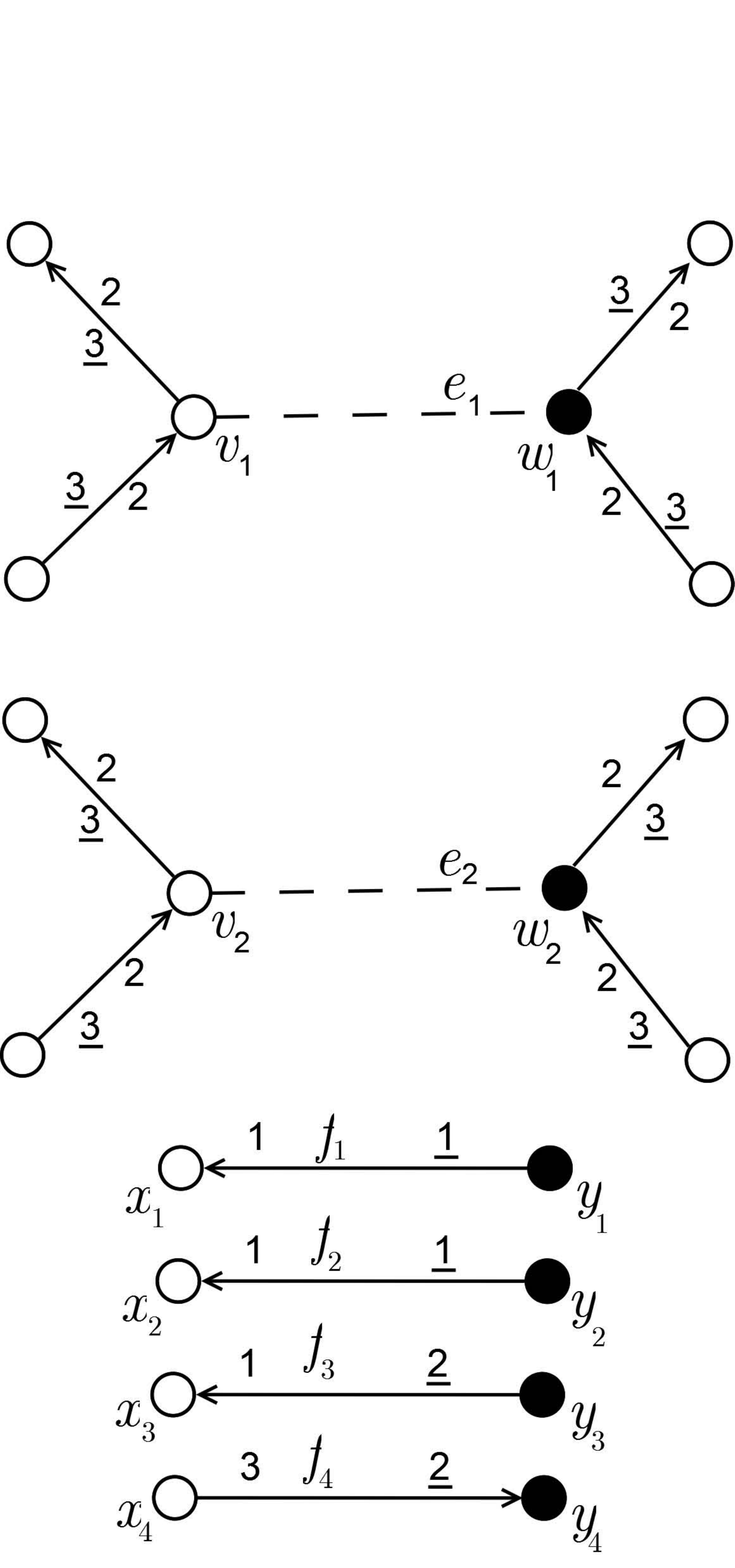}
 \caption{Cut}
   \label{Figure_Cut}
\end{figure}

Since $k=4$ and $S$ contains more white than black vertices, it follows 
with Lemma \ref{verschiedene2} that $v_1, v_2, x_1, x_2, x_3, x_4 \in A$. This implies 
that for $i \in \{1,2,3\}$ the edges $f_i$ are directed from $y_i$ to $x_i$ 
and $\phi(f_i) = 1$, and the edge $f_4$ is directed from $x_4$ to $y_4$, and  $\phi(f_i) = 3$.
The cut is depicted in Figure \ref{Figure_Cut}, where the underlined numbers are the colors and the other are the flow values of the edges.

\begin{claim} \label{S1}
Let $\{x,y\} = \{v_1,v_2\}$.
If $f_4 \not \in E(C_{x})$, then there is a directed path from $x$ to $y$ in $D(G[S])$.
\end{claim}
{\bf Proof.} Assume that $x=v_1$. 
Let $S'$ be the set of vertices $z$ in $D(G[S])$ such that there is a directed path from $v_1$ to $z$ in $D(G[S])$.
Suppose to the contrary that $v_2 \not \in S'$. 
Every edge of $\partial_{G[S]}(S')$ has its end vertex in $S'$. 
Since $\phi^-$ is a nowhere-zero flow on $G^-$ it follows, that $\delta\phi^-(S') = 0$.

Since $f_4$ and therefore $f_3$ as well are not in $E(C_{v_1})$
($f_1$ and and $f_2$ are colored with color 1 by $c'$), it follows that $V(C_{v_1}) \subseteq S'$.
But $f_4$ is the only outgoing edge of $\partial_{G^-}(S')$. 
Since $\phi^-(f_4) = 3$ and $\delta \phi^-(S') = 0$, it follows that $|\partial_{G^-}(S')| \leq 4$. Hence, 
$|\partial_{G}(S')| \leq 5$. Since each of $G[S']$ and $G[V(G)-S']$ contain odd circuits of a minimum 2-factor of $G$,
it follows that $m_G^* \leq 5$, a contradiction. $\Box$

Analogously we prove the following claim.

\begin{claim} \label{S2}
Let $\{x,y\} = \{w_1,w_2\}$.
If $f_4 \not \in E(C_{x})$, then there is a directed path from $y$ to $x$ in $D(G[V(G)-S])$.
\end{claim}
{\bf Proof.} Assume that $w_2=x$, and let $S'$ be the maximum set of vertices $z$ in $D(G[V(G)-S])$ 
such that there is a directed path from $z$ to $w_2$.
Suppose to the contrary that $w_1 \not \in S'$. Then all edges of $\partial_{G[V(G)-S]}(S')$ have their
initial end in $S'$. 
Since $f_4$ and therefore $f_3$ as well are not in $E(C_{w_2})$, it follows that $V(C_{w_2}) \subseteq S'$. 
But $f_4$ is the only ingoing edge of $\partial_{G^-}(S')$. 
Since $\phi^-(f_4) = 3$ and $\delta \phi^-(S') = 0$, it follows that $|\partial_{G^-}(S')| \leq 4$. Hence 
$|\partial_{G}(S')| \leq 5$. Since both $G[S']$ and $G[V(G)-S']$ contain odd circuits of a minimum 2-factor of $G$,
it follows that $m_G^* \leq 5$, a contradiction.
$\Box$

At most one circuit of $C_{v_1}$, $C_{v_2}$, $C_{w_1}$, $C_{w_2}$ contains the edge $f_4$ and therefore $f_3$ as well.
Suppose that $C_{v_2}$ contains $f_4$ (the other cases will be proved analogously). By Claim \ref{S1} there
is a directed path $P$ from $v_1$ to $v_2$ in $D(G[S])$. Since none of $C_{w_1}$ and $C_{w_2}$ contains $f_4$,
it follows by Claim \ref{S2} there is a directed path $P'$ from $w_2$ to $w_1$ in $D(G[V(G)-S])$. Hence $P$ and $P'$
are disjoint. Direct the edges $e_1$ and $e_2$ appropriately such that the
circuit $C$ with $E(C) = E(P) \cup E(P') \cup \{e_1,e_2\}$ is a directed circuit. Let $\phi_2$ be a 2-flow on $G$ with
$\phi_2(e)=1$, if $e \in E(C)$, and $\phi_2(e) = 0$ otherwise. Then $\phi^- + \phi_2$ is a nowhere-zero 5-flow on $G$, contradicting our 
supposition that $G$ has no nowhere-zero 5-flow. 

If none of the four odd circuits contains $f_4$, then the statement follows analogously. \hfill $\Box$

\end{document}